\documentclass[reqno]{amsart}

\usepackage{amssymb,amsmath,amsthm,amsfonts,fancyhdr, graphicx}
\usepackage{mathrsfs}
\usepackage{setspace}
\usepackage{cite}

\newtheorem{thm}{Theorem}[section]

\numberwithin{equation}{section}
\begin{document}

\title[exterior solutions to Monge-Amp\`{e}re equation]{A Remark on the asymptotic behavior of the exterior solutions to Monge-Amp\`{e}re equation}

\author{Guanghao Hong}
\address{School of Mathematics and Statistics, Xi'an Jiaotong University, Xi'an, P.R.China 710049.}
\email{ghhongmath@xjtu.edu.cn}

\begin{abstract}
We improve the result of Caffarelli-Li [CL03] on the asymptotic behavior at infinity of the exterior solution $u$ to Monge-Amp\`{e}re equation $det(D^2u)=1$ on $\mathbb{R}^n\backslash K$ for $n\geq 3$. We prove that the error term $O(|x|^{2-n})$ can be refined to $d (\sqrt{x'Ax})^{2-n}+O(|x|^{1-n})$ with $d=Res[u]$ the residue of $u$.
\end{abstract}

%\subjclass[2010]{Primary ; Secondary }


\maketitle

\section{Introduction}
The seminal results of J\"{o}rgens ($n=2$ [Jo54]), Calabi ($n\geq 5$ [Ca58]), and Pogorelov (all $n$ [Po72]) state that the classical convex solution to
\begin{equation*}
det(D^2u)=1 \ \ \ \mbox{in}\ \mathbb{R}^n
\end{equation*}
must be a quadratic polynomial. Caffarelli extended the result for classical solutions to viscosity solutions (see [CL03]).

Let $\mathcal{A}=\{A: A$ is $n\times n$ symmetric positive definite matrix with $det(A)=1\}$. Caffarelli and Yanyan Li [CL03] proved the following results.

\begin{thm}
Let $f\in C^0(\mathbb{R}^n)$ satisfies $0<\inf_{\mathbb{R}^n}f\leq \sup_{\mathbb{R}^n}f<\infty$ and the support of $f-1$ is bounded. Assume that $u$ is a convex viscosity solution of
\begin{equation*}
det(D^2u)=f \ \ \ \mbox{in}\ \mathbb{R}^n.
\end{equation*}
Then $u$ is $C^{\infty}$ in the complement of the support of $(f-1)$ and there exist some $A\in \mathcal{A}$, $b\in \mathbb{R}^n$ and $c\in \mathbb{R}$, such that

(i) for $n\geq 3$,
\begin{equation}
u(x)=\frac{1}{2}x'Ax+b\cdot x+c+O_k(|x|^{2-n}) \ \ \mbox{as}\ x\rightarrow \infty;
\end{equation}

(ii) for $n=2$,
\begin{equation}
u(x)=\frac{1}{2}x'Ax+b\cdot x+d\log \sqrt{x'Ax}+c+O_k(|x|^{-1}) \ \ \mbox{as}\ x\rightarrow \infty
\end{equation}
with $d=\frac{1}{2\pi}\int_{\mathbb{R}^2}(f-1)$.
The notation
$\varphi(x)=O_{k}(|x|^{m})$ means that $|D^{k}\varphi(x)|=O(|x|^{m-k})$ for all $k=0,1,2,\cdots$.
\end{thm}

\begin{thm}
Let $K$ be a bounded closed convex subset of $\mathbb{R}^n$, and let $u\in C^0(\mathbb{R}^n\backslash K)$ be a locally convex viscosity solution of
\begin{equation*}
det(D^2u)=1 \ \ \ \mbox{in}\ \mathbb{R}^n\backslash K
\end{equation*}
Then $u\in C^{\infty}(\mathbb{R}^n\backslash K)$ and there exist some $A\in \mathcal{A}$, $b\in \mathbb{R}^n$ and $c\in \mathbb{R}$, such that

(i) for $n\geq 3$, (1.1) holds;

(ii) for $n=2$, (1.2) holds for some $d\in \mathbb{R}$.
\end{thm}

We denote $u_i=\frac{\partial u}{\partial x_i}$, $u_{ij}=\frac{\partial^2 u}{\partial x_i\partial x_i}$ and $\tilde{u}_{ij}$ is the cofactor of $u_{ij}$. It is well known that the Monge-Amp\`{e}re operator has divergence structure
\begin{equation*}
det(D^2u)=\sum_{j=1}^n \partial_j(u_1\tilde{u}_{1j}):=div (\psi(u))
\end{equation*}
since the vector field $(\tilde{u}_{11},\tilde{u}_{12},\cdots,\tilde{u}_{1n})$ is divergence free (see \textit{e.g.} [BNST08]). Let $\xi(x)$ be any vector field in $\mathbb{R}^n$ satisfying $div \xi=1$, say $\xi(x)=x_1e_1$ or $\frac{x}{n}$. Let $U$ be a bounded domain with smooth boundary satisfying $U\supset supp (f-1)$ (in case of Theorem 1.1) or $U\supset K$ (in case of Theorem 1.2). Then the integral
$$\int_{\partial U} (\psi(u)-\xi)\cdot \vec{n}d\sigma=\int_{\partial U} \psi(u)\cdot \vec{n}d\sigma-|U|$$
 is independent of the specific choice of $U$. In case of Theorem 1.1, it is $\int_{\mathbb{R}^n}(f-1)$. We define
 $$Res[u]=\frac{1}{2\pi}\int_{\partial U} (\psi(u)-\xi)\cdot \vec{n}d\sigma \ \ \ \mbox{for}\ n=2,$$
  and
  $$Res[u]=\frac{1}{(n-2)n\omega_n}\int_{\partial U} (\psi(u)-\xi)\cdot \vec{n}d\sigma \ \ \ \mbox{for}\ n\geq 3,$$
   where $\omega_n$ denotes the volume of the unit ball in $\mathbb{R}^n$.

Therefore, in our notation, $d=Res[u]$ in (ii) of Theorem 1.1. Using the same method (p. 570 in [CL03]), one can also confirm that $d=Res[u]$ in (ii) of Theorem 1.2.

The residue $Res[u]$ is an essential quantity for $u$, so it is natural to expect that it also appears in the expansion of $u$ at infinity for $n\geq 3$. The purpose of this paper is to prove the following refined version of (1.1).

\begin{thm}
Under the conditions of Theorem 1.1 or Theorem 1.2, for $n\geq 3$, we have for some $A\in \mathcal{A}$, $b\in \mathbb{R}^n$ and $c\in \mathbb{R}$
\begin{equation*}
u(x)=\frac{1}{2}x'Ax+b\cdot x+c-Res[u](\sqrt{x'Ax})^{2-n}+O_k(|x|^{1-n}) \ \ \mbox{as}\ x\rightarrow \infty.
\end{equation*}
\end{thm}

\section{Proof of Theorem 1.3}
We prove Theorem 1.3 by an argument that we have used in [HY20] (see Step 3 in \S6.1).
\begin{proof}
Without loss of generality, we assume $A=I$, $b=0$ and $c=0$ because otherwise we can make some affine transformation as in [CL03]. Denote $$E(x):=u(x)-\frac{1}{2}|x|^2.$$ By (1.1), we have
\begin{equation}
E(x)=O(|x|^{2-n}),\ \ DE(x)=O(|x|^{1-n}),\ \ D^2E(x)=O(|x|^{-n}).
\end{equation}
We use the notation $F(\xi):=det(\xi_{ij})^{1/n}$. Then $E(x)$ satisfies the equation
\begin{equation}
\sum_{ij}a_{ij}(x)D_{ij}E(x)=F(I+D^2E)-F(I)=0
\end{equation}
in $\mathbb{R}^n\backslash B_{R_0}$ for some $R_0>0$, where $$a_{ij}(x)=\int_0^1 F_{\xi_{ij}}(I+sD^2E(x))ds=\delta_{ij}+O(|x|^{-n}).$$
We write (2.2) as
\begin{equation}
\triangle E(x)=\sum_{ij}(\delta_{ij}-a_{ij}(x))D_{ij}E(x):=g(x)=O(|x|^{-2n}).
\end{equation}

We use Kelvin transformation. Define $K[E](x):=|x|^{2-n}E(\frac{x}{|x|^2})$ for $x\in B_{\frac{1}{R_0}}\backslash \{0\}$. Then
\begin{equation}
\triangle K[E]=|x|^{-2-n}g(\frac{x}{|x|^2}):=\tilde{g}(x), \ \ \mbox{in}\ B_{\frac{1}{R_0}}\backslash \{0\}.
\end{equation}
From (2.3) and (2.4),we see $\tilde{g}(x)=O(|x|^{n-2})$ as $x\rightarrow 0$, so $\tilde{g}(x)\in L^{\infty}(B_{\frac{1}{R_0}})$. Let $N[\tilde{g}]$ be the  Newtonian potential of $\tilde{g}$
in $B_{\frac{1}{R_{0}}}$. Since $\tilde{g}$ is in $L^{p}(B_{\frac{1}{R_{0}}})$ for
any $p>0$, $N[\tilde{g}]$ is $W^{2,p}$ for any $p$ and hence is $C^{1,\alpha}$
for any $0<\alpha<1$. Now $K[E]-N[\tilde{g}]$ is harmonic in $B_{\frac{1}{R_{0}}}\backslash\{0\}$. From (2.1), we know $K[E]$ is bounded. So $K[E]-N[\tilde{g}]$ is bounded and hence $\{0\}$ is its removable singularity. That is $K[E]-N[\tilde{g}]$ is harmonic in $B_{\frac{1}{R_{0}}}$. So $K[E](x)$ is a $C^{1,\alpha}$
function in $B_{\frac{1}{R_{0}}}$. Fix an $\alpha\in(0,1)$, for some affine
function $\tilde{c}+\tilde{b}\cdot x$, we have $$|K[E](x)-\tilde{c}-\tilde{b}\cdot x|\leq
C|x|^{1+\alpha}\ \ \ \mbox{in}\ B_{\frac{1}{R_{0}}}.$$ Going back to $E$,
we have $$|E(x)-\tilde{c}|x|^{2-n}-\tilde{b}\cdot\frac{x}{|x|^n}|\leq
C|x|^{1-n-\alpha}\ \ \ \mbox{for}\ |x|\geq R_{0}.$$
That is
\begin{equation*}
E(x)=\tilde{c}|x|^{2-n}+O(|x|^{1-n}) \ \ \ \mbox{as}\ x\rightarrow \infty.
\end{equation*}

By Lemma 3.5 in [CL03], one can improve the above $O(|x|^{1-n})$ to $O_k(|x|^{1-n})$. The remaining thing is to confirm that $\tilde{c}=-Res[u]$. This can be done in the same way as in the 2 dimensional case (p.570 in [CL03]). We give the details in the following.

For $r>R_0$,
\begin{eqnarray}
Res[u]&=&\frac{1}{(n-2)n\omega_n}\int_{\partial B_r} (\psi(u)-x_1e_1)\cdot \vec{n}d\sigma\nonumber \\
      &=&\frac{1}{(n-2)n\omega_n r}\int_{\partial B_r} (\sum_{j=1}^n u_1 \tilde{u}_{1j} x_j-x_1^2) d\sigma.
\end{eqnarray}
We write $u=w+\eta+O_k(|x|^{1-n})$, where $w=\frac{|x|^2}{2}$ and $\eta=\tilde{c}|x|^{2-n}$. By computation, we have
\begin{eqnarray*}
u_1&=&w_1+\eta_1+O(r^{-n}) \\
      &=&x_1+\tilde{c}(2-n)r^{-n}x_1+O(r^{-n}),
\end{eqnarray*}
\begin{eqnarray*}
\tilde{u}_{11}&=&1+\sum_{k=2}^n \eta_{kk}+O(r^{-2n}) \\
      &=&1-\eta_{11}+O(r^{-2n})\ \ \ \ \ \ (\mbox{since}\ \triangle \eta=0)\\
      &=&1-\tilde{c}(2-n)r^{-n}-\tilde{c}(2-n)(-n)r^{-n-2}x_1^2+O(r^{-2n})
\end{eqnarray*}
and for $j\geq 2$,
\begin{eqnarray*}
\tilde{u}_{1j}&=&-\eta_{j1}+O(r^{-2n}) \\
      &=&-\tilde{c}(2-n)(-n)r^{-n-2}x_1x_j+O(r^{-2n}).
\end{eqnarray*}
So
\begin{eqnarray}
\sum_{j=1}^n u_1 \tilde{u}_{1j} x_j&=&x_1^2+x_1\eta_1-\eta_{11}x_1^2-\sum_{j=2}^n x_1\eta_{j1}x_j+O(r^{1-n})\nonumber \\
&=&x_1^2-\tilde{c}n(n-2)r^{-n}x_1^2+O(r^{1-n}).
\end{eqnarray}
Inserting (2.6) into (2.5), we have
\begin{eqnarray*}
Res[u]=-\tilde{c}+O(r^{-1}).
\end{eqnarray*}
Letting $r\rightarrow\infty$, we get $$Res[u]=-\tilde{c}.$$
\end{proof}

\end{document}